\documentstyle[txmac,a4,
amssymb,%
case,%
twoside,%
nocaphead,%
epsf,%
mypic,times,mathptm]{article}

\advance\oddsidemargin by -1.9cm
\advance\evensidemargin by -1.9cm
\advance\textwidth by 3.8cm

\def\mynewtheo#1#2{%
\newtheorem{@#1}{#2}[section]%
\newenvironment{#1}{\begin{@#1}\rm}{\end{@#1}}}

\mynewtheo{lemma}{Lemma}
\mynewtheo{exer}{Exercise}
\mynewtheo{theo}{Theorem}
\mynewtheo{rem}{Remark}
\mynewtheo{defi}{Definition}
\mynewtheo{conj}{Conjecture}
\mynewtheo{corr}{Corollary}
\mynewtheo{prop}{Proposition}
\mynewtheo{question}{Question}
\mynewtheo{exam}{Example}

\newenvironment{theorem}{\begin{theo}}{\end{theo}}

\newenvironment{eqn}{\begin{equation}}{\end{equation}}

\begin{document}

\makeatletter

\newenvironment{myeqn*}[1]{\begingroup\def\@eqnnum{\reset@font\rm#1}%
\xdef\@tempk{\arabic{equation}}\begin{equation}\edef\@currentlabel{#1}}
{\end{equation}\endgroup\setcounter{equation}{\@tempk}\ignorespaces}

\newenvironment{myeqn}[1]{\begingroup\let\eq@num\@eqnnum
\def\@eqnnum{\bgroup\let\r@fn\normalcolor 
\def\normalcolor####1(####2){\r@fn####1#1}%
\eq@num\egroup}%
\xdef\@tempk{\arabic{equation}}\begin{equation}\edef\@currentlabel{#1}}
{\end{equation}\endgroup\setcounter{equation}{\@tempk}\ignorespaces}

\newcommand{\mybin}[2]{\text{$\Bigl(\begin{array}{@{}c@{}}#1\\#2%
\end{array}\Bigr)$}}
\newcommand{\mybinn}[2]{\text{$\biggl(\begin{array}{@{}c@{}}%
#1\\#2\end{array}\biggr)$}}

\def\overtwo#1{\mbox{\small$\mybin{#1}{2}$}}
\newcommand{\mybr}[2]{\text{$\Bigl\lfloor\mbox{%
\small$\displaystyle\frac{#1}{#2}$}\Bigr\rfloor$}}
\def\mybrtwo#1{\mbox{\mybr{#1}{2}}}

\def\myfrac#1#2{\raisebox{0.2em}{\small$#1$}\!/\!\raisebox{-0.2em}{\small$#2$}}

\def\myeqnlabel{\bgroup\@ifnextchar[{\@maketheeq}{\immediate
\stepcounter{equation}\@myeqnlabel}}

\def\@maketheeq[#1]{\def\theequation{#1}\@myeqnlabel}

\def\@myeqnlabel#1{%
{\edef\@currentlabel{\theequation}
\label{#1}\enspace\eqref{#1}}\egroup}

\def\rato#1{\hbox to #1{\rightarrowfill}}
\def\arrowname#1{{\enspace
\setbox7=\hbox{F}\setbox6=\hbox{%
\setbox0=\hbox{\footnotesize $#1$}\setbox1=\hbox{$\to$}%
\dimen@\wd0\advance\dimen@ by 0.66\wd1\relax
$\stackrel{\rato{\dimen@}}{\copy0}$}%
\ifdim\ht6>\ht7\dimen@\ht7\advance\dimen@ by -\ht6\else
\dimen@\z@\fi\raise\dimen@\box6\enspace}}

\def\ptat#1{{\picfillgraycol{0}\picfilledcircle{#1}{0.05}{}}}
\def\chrd#1#2{\picline{1 #1 polar}{1 #2 polar}}
\def\arrow#1#2{\picvecline{1 #1 polar}{1 #2 polar}}

\def\labch#1#2#3{\chrd{#1}{#2}\picputtext{1.3 #2 polar}{$#3$}}
\def\labar#1#2#3{\arrow{#1}{#2}\picputtext{1.3 #2 polar}{$#3$}}
\def\labbr#1#2#3{\arrow{#1}{#2}\picputtext{1.3 #1 polar}{$#3$}}

\def\CD{\szCD{4mm}}
\def\szCD#1#2{{\let\@nomath\@gobble\small\diag{#1}{2.4}{2.4}{
  \pictranslate{1.2 1.2}{
    \piccircle{0 0}{1}{}
    #2
}}}}


\author{A. Stoimenow\footnotemark[1]\\[2mm]
\small Ludwig-Maximilians University Munich, Mathematics\\
\small Institute, Theresienstra\ss e 39, 80333 M\"unchen, Germany,\\
\small e-mail: {\tt stoimeno@informatik.hu-berlin.de},\\
\small WWW\footnotemark[2]\enspace: {\hbox{\tt http://www.informatik.hu-berlin.de/%
\raisebox{-0.8ex}{\tt\~{}}stoimeno}}
}

{\def\thefootnote{\fnsymbol{footnote}}
\footnotetext[1]{Supported by a DFG postdoc grant.}
\footnotetext[2]{All my papers, including those referenced here, are
available on my webpage or by sending me an inquiry to the email
address I specified above. Hence, \em{please}, do not complain about
some paper being non-available before trying these two options. Thank
you.}
}

\title{\large\bf \uppercase{The crossing number and maximal bridge
length of a knot diagram}\\[4mm]
{\it\small This is a preprint. I would be grateful
for any comments and corrections!}}

\date{\large Current version: \today\ \ \ First version:
\makedate{15}{2}{1999}}

\maketitle

\makeatletter

\def\pt#1{{\picfillgraycol{0}\picfilledcircle{#1}{0.06}{}}}
\def\labpt#1#2#3{\pictranslate{#1}{\pt{0 0}\picputtext{#2}{$#3$}}}

\let\vn\varnothing
\let\point\pt
\let\ay\asymp
\let\pa\partial
\let\al\alpha
\let\be\beta
\let\Dl\Delta
\let\Gm\Gamma
\let\gm\gamma
\let\de\delta
\let\dl\delta
\let\eps\epsilon
\let\lm\lambda
\let\Lm\Lambda
\let\sg\sigma
\let\vp\varphi
\let\om\omega
\let\diagram\diag

\let\sm\setminus
\let\tl\tilde
\def\ncap{\not\mathrel{\cap}}
\def\sgn{\text{\rm sgn}\,}
\def\cf{\text{\rm cf}\,}
\def\md{\max\deg}
\def\mc{\max\cf}
\def\Lra{\Longrightarrow}
\def\lra{\longrightarrow}
\def\so{\Rightarrow}
\def\So{\Longrightarrow}
\def\nin{\not\in}
\let\ds\displaystyle
\def\bt{\bar t_2}
\def\is{\big|\,S\cap [-k,k]^{n}\,\big|}
\let\llra\longleftrightarrow
\let\reference\ref

\long\def\@makecaption#1#2{%
   \vskip 10pt
   {\let\label\@gobble
   \let\ignorespaces\@empty
   \xdef\@tempt{#2}%
   }%
   \ea\@ifempty\ea{\@tempt}{%
   \setbox\@tempboxa\hbox{%
      \fignr#1#2}%
      }{%
   \setbox\@tempboxa\hbox{%
      {\fignr#1:}\capt\ #2}%
      }%
   \ifdim \wd\@tempboxa >\captionwidth {%
      \rightskip=\@captionmargin\leftskip=\@captionmargin
      \unhbox\@tempboxa\par}%
   \else
      \hbox to\captionwidth{\hfil\box\@tempboxa\hfil}%
   \fi}%
\def\fignr{\small\sffamily\bfseries}%
\def\capt{\small\sffamily}%

\newdimen\@captionmargin\@captionmargin2cm\relax
\newdimen\captionwidth\captionwidth\hsize\relax

\def\eqref#1{(\protect\ref{#1})}

\def\proof{\@ifnextchar[{\@proof}{\@proof[\unskip]}}
\def\@proof[#1]{\noindent{\bf Proof #1.}\enspace}

\def\hint{\noindent Hint: }
\def\problem{\noindent{\bf Problem.} }

\def\@mt#1{\ifmmode#1\else$#1$\fi}
\def\qed{\hfill\@mt{\Box}}
\def\qqed{\hfill\@mt{\Box\enspace\Box}}

\def\cU{{\cal U}}
\def\cC{{\cal C}}
\def\cP{{\cal P}}
\def\tP{{\tilde P}}
\def\tZ{{\tilde Z}}
\def\fg{{\frak g}}
\def\tr{\text{tr}}
\def\cZ{{\cal Z}}
\def\cD{{\cal D}}
\def\bR{{\Bbb R}}
\def\cE{{\cal E}}
\def\bZ{{\Bbb Z}}
\def\bN{{\Bbb N}}

\def\bysame{\same[\kern2cm]\,}

\def\br#1{\left\lfloor#1\right\rfloor}
\def\BR#1{\left\lceil#1\right\rceil}

\def\abstractname{}

\@addtoreset {footnote}{page}

\renewcommand{\section}{%
   \@startsection
         {section}{1}{\z@}{-1.5ex \@plus -1ex \@minus -.2ex}%
               {1ex \@plus.2ex}{\large\bf}%
}
\renewcommand{\@seccntformat}[1]{\csname the#1\endcsname .
\quad}

\def\bC{{\Bbb C}}
\def\bP{{\Bbb P}}

\def\@test#1#2#3#4{%
  \let\@tempa\go@
  \@tempdima#1\relax\@tempdimb#3\@tempdima\relax\@tempdima#4\unitxsize\relax
  \ifdim \@tempdimb>\z@\relax
    \ifdim \@tempdimb<#2%
      \def\@tempa{\@test{#1}{#2}}%
    \fi
  \fi
  \@tempa
}

\def\go@#1\@end{}
\newdimen\unitxsize
\newif\ifautoepsf\autoepsftrue

\unitxsize4cm\relax
\def\epsfsize#1#2{\epsfxsize\relax\ifautoepsf
  {\@test{#1}{#2}{0.1 }{4   }
		{0.2 }{3   }
		{0.3 }{2   }
		{0.4 }{1.7 }
		{0.5 }{1.5 }
		{0.6 }{1.4 }
		{0.7 }{1.3 }
		{0.8 }{1.2 }
		{0.9 }{1.1 }
		{1.1 }{1.  }
		{1.2 }{0.9 }
		{1.4 }{0.8 }
		{1.6 }{0.75}
		{2.  }{0.7 }
		{2.25}{0.6 }
		{3   }{0.55}
		{5   }{0.5 }
		{10  }{0.33}
		{-1  }{0.25}\@end
		\ea}\ea\epsfxsize\the\@tempdima\relax
		\fi
		}

\def\rr#1{\unitxsize2.8cm\relax\epsffile{#1.eps}}
\def\crr#1{\rr{curve#1}}
\def\vis#1#2{\hbox{\begin{tabular}{c}\rr{t-#1-#2} \\ $#1_{#2}$\end{tabular}}}


\makeatletter
\def\mybrace#1#2{\@tempdima#1em\relax
\advance\@tempdima by -1em\relax
\setbox\@tempboxa=\hbox{\raisebox{-0.5\@tempdima}{$\ds 
\left.\rule[0.5\@tempdima]{\z@
}{0.5\@tempdima}\right\} #2$}}\dp\@tempboxa=\z@
\box\@tempboxa}

\def\namedarrow#1{{\enspace
\setbox7=\hbox{F}\setbox6=\hbox{%
\setbox0=\hbox{\footnotesize $#1$}\setbox1=\hbox{$\to$}%
\dimen@\wd0\advance\dimen@ by 0.66\wd1\relax
$\stackrel{\rato{\dimen@}}{\copy0}$}%
\ifdim\ht6>\ht7\dimen@\ht7\advance\dimen@ by -\ht6\else
\dimen@\z@\fi\raise\dimen@\box6\enspace}}

\let\old@tl\~\def\~{\raisebox{-0.8ex}{\tt\old@tl{}}}
\let\lra\longrightarrow
\def\bt{\bar t'_{2}}
\let\sm\setminus
\let\eps\varepsilon
\let\ex\exists
\let\fa\forall
\let\ps\supset

\def\rs#1{\raisebox{-0.4em}{$\big|_{#1}$}}


\def\@test#1#2#3#4{%
  \let\@tempa\go@
  \@tempdima#1\relax\@tempdimb#3\@tempdima\relax\@tempdima#4\unitxsize\relax
  \ifdim \@tempdimb>\z@\relax
    \ifdim \@tempdimb<#2%
      \def\@tempa{\@test{#1}{#2}}%
    \fi
  \fi
  \@tempa
}

\newbox\@tempboxb
\def\epsfs{\@ifnextchar[{\@epsfs}{\@@epsfs}}
\def\@epsfs[#1]#2{{\ifautoepsf\unitxsize#1\relax\else
\epsfxsize#1\relax\fi\@@epsfs{#2}}}
\def\@@epsfs#1{\setbox\@tempboxb=\hbox{\,\epsffile{#1.eps}}\,
\parbox{\wd\@tempboxb}{\box\@tempboxb}}
	
\def\eepsfs{\@ifnextchar[{\@eepsfs}{\@@eepsfs}}
\def\@eepsfs[#1]#2{\uu{\ifautoepsf\unitxsize#1\relax\else
\epsfxsize#1\relax\fi\epsffile{#2.eps}}}
\def\@@eepsfs#1{\uu{\epsffile{#1.eps}}}

{\let\@noitemerr\relax
\vskip-2.7em\kern0pt\begin{abstract}
\noindent{\bf Abstract.}\enspace
We give examples showing that Kidwell's inequality for the
maximal degree of the Brandt-Lickorish-Millett-Ho polynomial is
in general not sharp.
\end{abstract}
}\vspace{7mm}

{\parskip0.2mm\tableofcontents}
\vspace{7mm}

\section{Introduction}

\begin{figure}
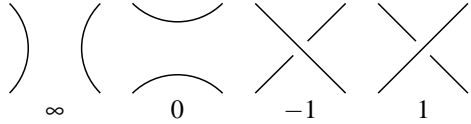

\[
\begin{array}{*4c}
\diag{6mm}{2}{2}{
  \piccirclearc{-1 1}{1.41}{-45 45}
  \piccirclearc{3 1}{1.41}{135 -135}
} & 
\diag{6mm}{2}{2}{
  \piccirclearc{1 -1}{1.41}{45 135}
  \piccirclearc{1 3}{1.41}{-135 -45}
} &
\diag{6mm}{2}{2}{
\picmultiline{-13.5 1 -1 0}{0 0}{2 2}
\picmultiline{-13.5 1 -1 0}{2 0}{0 2}
} &
\diag{6mm}{2}{2}{
\picmultiline{-13.5 1 -1 0}{2 0}{0 2}
\picmultiline{-13.5 1 -1 0}{0 0}{2 2}
}
\\[2mm]
\infty & 0 & -1 & 1 
\end{array}
\]
\caption{\label{figtan}The Conway tangles.}
\end{figure}

The $Q$ (or absolute) polynomial is a polynomial invariant in one
variable $z$ of links
(and in particular knots) in $S^3$ without orientation defined
by being 1 on the unknot and the relation
\begin{eqn}\label{Qrel}
A_{1}+A_{-1}=z(A_{0}+A_{\infty})\,,
\end{eqn}
where $A_i$ are the $Q$ polynomials of links $K_i$ and $K_i$ ($i\in\bZ
\cup\{\infty\}$) possess diagrams equal except in one room, where
an $i$-tangle (in the Conway \cite{Conway} sense) is inserted,
see figure \reference{figtan}.

It has been discovered in 1985 independently by Brandt, Lickorish
and Millett \cite{BLM} and Ho \cite{Ho}. Several months after
its discovery, Kauffman \cite{Kauffman} found a 2-variable polynomial
$F(a,z)$, specializing to $Q$ by setting $a=1$.

In \cite{Kidwell}, Kidwell found a nice inequality for the maximal
degree of the $Q$ polynomial.

\begin{theorem} (Kidwell)
Let $D$ be a diagram of a knot (or link) $K$. Then 
\begin{eqn}\label{Qineq}
\md Q(K)\,\le\,c(D)-d(D)\,,
\end{eqn}
where $c(D)$ is the crossing number of $D$ and $d(D)$ its maximal
bridge length, i.~e., the maximal number of consecutive crossing
over- or underpasses. Moreover, if $D$ is alternating (i.~e. $d(D)=1$)
and prime, then equality holds in \eqref{Qineq}.
\end{theorem}

In \cite[problem 4, p.~560]{Morton} he asked whether \eqref{Qineq}
always becomes equality when minimizing the r.h.s. over all diagrams $D$
of $K$. From the theorem it follows that this is true for alternating
knots and also for those non-alternating knots $K$, where
$\md Q(K)=c(K)-2$ (here $c(K)$ denotes the crossing number of $K$).
All non-alternating knots in Rolfsen's tables \cite{Rolfsen}
have this property except for one~-- the Perko knot $10_{161}$
(and its obversed duplication $10_{162}$), where $\md Q=6$. Hence,
as quoted by Kidwell, this knot became a promising candidate for
strict inequality in \eqref{Qineq}. To express ourselves more easily,
we define

\begin{defi}
Call a knot $K$ $Q$-maximal, if \eqref{Qineq} with the r.h.s. minimized
over all diagrams $D$ of $K$ becomes equality.
\end{defi}

The aim of this note is to show that indeed the Perko knot is not
$Q$-maximal. We give several modifications of our arguments and
examples showing how they can be applied to exhibit non-$Q$-maximality.

\section{\label{S2}Plane curves}

We start by some discussion on plane curves.

\begin{defi}
A non-closed plane curve is a $C^1$ map $\gm\,:\,[0,1]\,\to\,\bR^2$
with $\gm(0)\ne \gm(1)$ and only transverse self-intersections. $\gm$
carries a natural orientation.
\end{defi}

\begin{exam} Here are some plane curves:
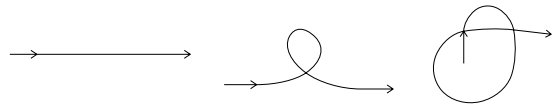
\[
\epsfs[6mm]{t-curve0}\quad
\epsfs[1.5cm]{t-curve1}\quad
\epsfs[1.5cm]{t-curve2}
\]
\end{exam}

In the following, whenever talking of plane curves we mean non-closed
ones with orientation unless otherwise stated. However, in some
cases it is possible to forget about orientation if it is irrelevant.
It is convenient to
identify $\gm$ with $\gm([0,1])$ wherever this causes no confusion.
Whenever we want to emphasize that a line segment in a local picture
starts with an endpoint, the endpoint will be depicted as a thickened
dot.

\begin{defi}
The crossing number $c(\gm)$ of a curve $\gm$ is the number of
self-intersections (crossings). The curve $\gm$ with $c(\gm)=0$
is called trivial.
\end{defi}

\begin{defi}
We call a non-closed curve $\tl\gm$ similar to $\gm$ (and denote it
by $\tl\gm\sim\gm$) if $\tl\gm(0)=\gm(0)$, $\tl\gm(1)=\gm(1)$ and
$\tl\gm$ intersects $\gm$ only transversely. The distance $d(\gm)$ of $\gm$
call the number $\min\{\,\#(\gm\cap\tl\gm)\,:\,\tl\gm\sim\gm\,\}-2$ (the
`$-2$' provided to ignore the coincidence of start- end endpoint).
A curve $\tl\gm$ realizing this minimum is called minimal similar
to $\gm$. Such $\tl\gm$ can be chosen to have no self-intersections.
\end{defi}

\begin{exam}
The curves $\epsfs[2mm]{t-curve0}$ and
$\epsfs[0.6cm]{t-curve1}$ have $d=0$, while $d\left(\epsfs[0.6cm]{t-curve2}
\right)=1$.
\end{exam}

\begin{defi}
We call a plane curve $\gm$ composite, if there is a \em{closed}
plane curve $\gm'$ (with no self-intersections) such that $\gm'$
intersects $\gm$ in exactly one point, transversely, and in both
components of $\bR^2\sm\gm'$ there are crossings of $\gm$. In this case
$\gm'$ separates $\gm$ into two parts $\gm_1$ and $\gm_2$, which we
call components of $\gm$. We write $\gm=\gm_1\#\gm_2$. Conversely,
this can be used to define the operation `$\#$' (connected sum)
of $\gm_1$ and $\gm_2$, wherever $\gm_1(1)$ or $\gm_2(0)$ are
in the closure of the unbounded component of their complements.
We call $\gm$ prime, if it is not composite.
\end{defi}

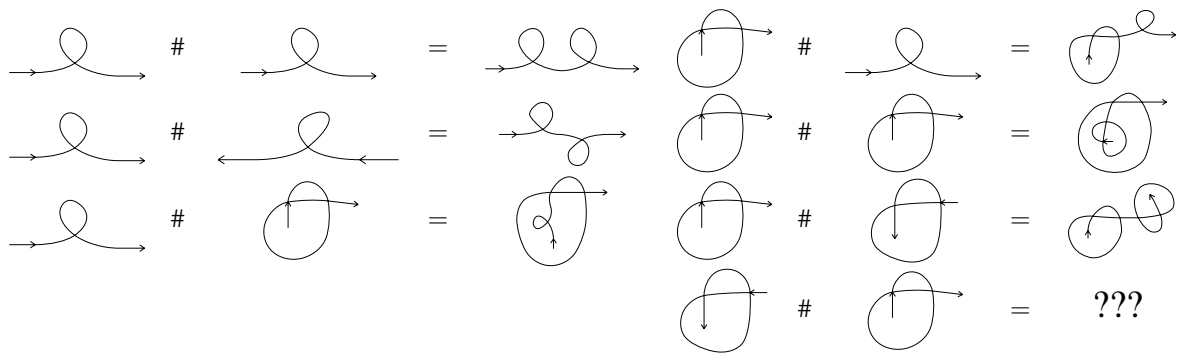
\begin{exam}
\[
\begin{array}{*{10}c}
\epsfs[1.2cm]{t-curve1}&\#&\epsfs[1.2cm]{t-curve1} & \,=\, &
\epsfs[1.2cm]{t-curve1.1} &
\epsfs[1.2cm]{t-curve2}&\#&\epsfs[1.2cm]{t-curve1} & \,=\, &
\epsfs[1.2cm]{t-curve2.1} \\
\epsfs[1.2cm]{t-curve1}&\#&\epsfs[1.2cm]{t-curve1-} & \,=\, &
\epsfs[1.2cm]{t-curve1.1-} &
\epsfs[1.2cm]{t-curve2}&\#&\epsfs[1.2cm]{t-curve2} & \,=\, &
\epsfs[1.2cm]{t-curve2.2} \\
\epsfs[1.2cm]{t-curve1}&\#&\epsfs[1.2cm]{t-curve2} & \,=\, &
\epsfs[1.2cm]{t-curve1.2} &
\epsfs[1.2cm]{t-curve2}&\#&\epsfs[1.2cm]{t-curve2-} & \,=\, &
\epsfs[1.2cm]{t-curve2.2-} \\
& & & & & 
\epsfs[1.2cm]{t-curve2-}&\#&\epsfs[1.2cm]{t-curve2} & \,=\, &
\mbox{\Large ???}
\end{array}
\]
\end{exam}
This example visualizes that the connected sum in general depends on
the orientation of the summands and their order.

It is clear that the crossing number is additive under connected sum and
it's a little exercise to verify that the distance is as well.

The path we are going to follow starts with the following

\begin{exer}\label{x1}
Verify that the complement of a curve $\gm$ has $c(\gm)+1$
connected components, and conclude from that that $d(\gm)\le c(\gm)$.
\end{exer}

\hint One way to show that is to observe it when $\gm$ is trivial,
to prove that you can obtain any $\gm$ from the trivial one
by the four local moves
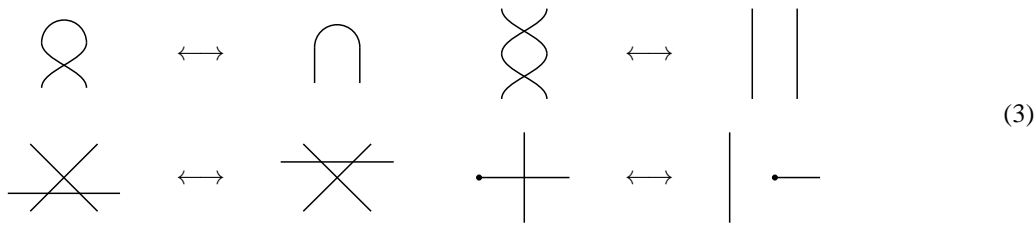
\begin{equation}\label{Reimov}
\begin{array}{*6c}
\diag{6mm}{1.5}{1.5}
{\pictranslate{0.25 0}{
   \braid{0.5 0.5}{1 1}
   }
\piccurve{0.25 1}{0.25 1.25}{0.45 1.5}{0.75 1.5}
\pictranslate{1.5 0}{\picscale{-1 1}{
  \piccurve{0.25 1}{0.25 1.25}{0.45 1.5}{0.75 1.5}
  }}
}\quad&\llra&\quad
\diag{6mm}{1.5}{1.3}
{\pictranslate{0.25 0}{
\picline{0 0}{0 0.8}
\picline{1 0}{1 0.8}
\piccirclearc{0.5 0.8}{0.5}{0 180}
}}&\qquad
\diag{6mm}{1}{2}{
\braid{0.5 0.5}{1 1}
\braid{0.5 1.5}{1 1}
}\quad&\llra&\quad
\diag{6mm}{1}{2}{
\picline{0 0}{0 2}
\picline{1 0}{1 2}
}
\\[8mm]
\diag{6mm}{2.5}{1.5}{
\picline{0 0.4}{2.5 0.4}
\picline{2 1.5}{0.5 0}
\picline{2 0}{0.5 1.5}
}\quad&\llra&\quad
\diag{6mm}{2.5}{1.5}{
\picline{0 1.1}{2.5 1.1}
\picline{2 1.5}{0.5 0}
\picline{2 0}{0.5 1.5}
}&\qquad
\diag{6mm}{2}{2}{
  \picline{1 2}{1 0}
  \picline{0 1}{2 1}
  \ptat{0 1}
}\quad&\llra&\quad
\diag{6mm}{2}{2}{
  \picline{0 2}{0 0}
  \picline{1 1}{2 1}
  \ptat{1 1}
}
\end{array}
\end{equation}
%
%
and to trace how the number of components and $c(\gm)$ change under
these moves.

Our first aim is to improve slightly the inequality of the exercise.

\begin{lemma}\label{l1}
$d(\gm)\le\max(2,c(\gm)-2)$ if $\gm$ prime.
\end{lemma}

\proof Consider the first (and analogously last) crossing of $\gm$
(that is, the crossings passed as first and last by $\gm$).
Denote by letters the connected components of the complement near
these crossings:
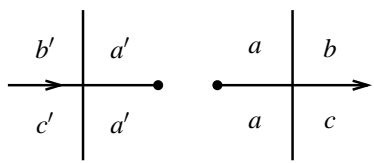
\[
\diag{1cm}{2}{2}{
  \picline{1 2}{1 0}
  \picvecline{0 1}{0.7 1}
  \picline{0.7 1}{2 1}
  \ptat{2 1}
  \picputtext{1.5 1.5}{$a'$}
  \picputtext{1.5 0.5}{$a'$}
  \picputtext{0.5 1.5}{$b'$}
  \picputtext{0.5 0.5}{$c'$}
}\qquad
\diag{1cm}{2}{2}{
  \picline{1 2}{1 0}
  \picvecline{0 1}{2 1}
  \ptat{0 1}
  \picputtext{1.5 1.5}{$b$}
  \picputtext{1.5 0.5}{$c$}
  \picputtext{0.5 1.5}{$a$}
  \picputtext{0.5 0.5}{$a$}
}\]
First note that $a\nin\{b,c\}$, else if w.l.o.g. $a=c$ there would
be a closed curve $\gm'$ like
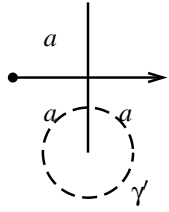
\[
\diag{1cm}{2}{3}{
  \pictranslate{0 1}{
    \picline{1 2}{1 0}
    \picvecline{0 1}{2 1}
    \ptat{0 1}
    \picputtext{1.5 0.5}{$a$}
    \picputtext{0.5 1.5}{$a$}
    \picputtext{0.5 0.5}{$a$}
  }%
  \piclinedash{0.2 0.1}{0.25}
  \piccircle{1 1}{0.6}{}
  \picputtext{1.7 0.5}{$\gm'$}
}
\]
intersecting $\gm$ in only one point 
and either $d(\gm)=0$ or $\gm$ is composite.

Then note that $b\ne c$, because else there would be a $\gm'$ like
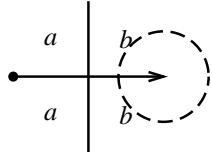
\[
\diag{1cm}{3}{2}{
  \picline{1 2}{1 0}
  \picvecline{0 1}{2 1}
  \ptat{0 1}
  \picputtext{1.5 1.5}{$b$}
  \picputtext{1.5 0.5}{$b$}
  \picputtext{0.5 1.5}{$a$}
  \picputtext{0.5 0.5}{$a$}
  \piclinedash{0.2 0.1}{0.25}
  \piccircle{2 1}{0.6}{}
}
\]
and the 2 curve segments could not be connected.  

Therefore, $b\ne c$ and a minimal curve $\tl\gm\sim\gm$ would not
need to pass through one of $b$ and $c$. The same holds for the
last crossing of $\gm$. Hence we avoid $\tl\gm$ passing through
at least two components of the complement of $\gm$, unless $\{b,c\}
\cap\{b',c'\}\ne\vn$, but then $d(\gm)\le 2$. \qed

For the Perko knot we need to work a little harder.

\begin{lemma}\label{l2}
$d(\gm)\le\max(3,c(\gm)-3)$ if $\gm$ prime.
\end{lemma}

\begin{defi}
An isolated crossing of $\gm$ is a crossing $p$ such that there
is a closed curve $\gm'$ with $\gm\cap\gm'=\{p\}$ and $\gm'$
intersects transversely both strands of $\gm$ intersecting at $p$.
\end{defi}

\proof[of lemma] If $\gm$ has an isolated crossing, then one
of the components of $\bR^2\sm\gm$ has both and the other one has
no one of the endpoints of $\gm$. Removing the part of $\gm$ in
latter component and smoothing $\gm$ near $p$ reduces $c(\gm)$,
but not $d(\gm)$, hence we may (say, by induction on $c(\gm)$)
assume that $\gm$ has no isolated crossing.

Now consider a crossing of $\gm$ which is neither the first nor
the last and denote the components near it by $l$, $m$, $n$ and $o$.
\begin{eqn}\label{crossing}
\diag{1cm}{2}{2}{
    \picline{1 2}{1 0}
    \picline{0 1}{2 1}
    \picputtext{1.5 1.5}{$m$}
    \picputtext{1.5 0.5}{$o$}
    \picputtext{0.5 1.5}{$l$}
    \picputtext{0.5 0.5}{$n$}
}
\end{eqn}
Call 2 components neighbored if the intersection of the
closures of their fragments in \eqref{crossing} is a line, and opposite
if it is just the crossing itself.

By primality of $\gm$ any two neighbored components are distinct and
by non-isolatedness of the crossing so are any two opposite components.

Hence $l$, $m$, $n$ and $o$ are pairwise distinct. Now call $b$, $b'$
the components which were found not to be passed by a minimal similar
curve $\gm'$ to $\gm$ by the proof of lemma \reference{l1} and $a$, $a'$
the components denoted so in the same proof.

Then distinguish some cases.

\begin{caselist}
\case\label{c1}
No one of $b$, $b'$ is among $l$, $m$, $n$ and $o$. As $\gm'$
certainly does not pass through one of $l$, $m$, $n$ and $o$ you
have a third component not passed by $\gm'$ and you are done.

\case Exactly one of $b$, $b'$, say $b$, is among $l$, $m$, $n$ and $o$.
You would be done as in case \reference{c1} unless $\gm'$ does
not pass only through $b$. Then you have a picture like this:
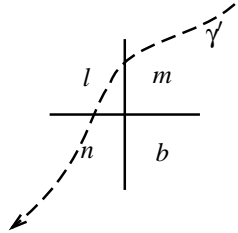
\[
\diag{1cm}{3}{3}{
  \pictranslate{0.5 0.5}{
    \picline{1 2}{1 0}
    \picline{0 1}{2 1}
    \picputtext{1.5 1.5}{$m$}
    \picputtext{1.5 0.5}{$b$}
    \picputtext{0.5 1.5}{$l$}
    \picputtext{0.5 0.5}{$n$}
  }
  \piclinedash{0.2 0.1}{0.25}
  \picstroke{
    \piccurve{3 3}{2.5 2.5}{1.5 2.5}{1.3 1.9}
    \picveccurveto{1 1.5}{1 0.8}{0 0}
  }
  \picputtext{2.7 2.6}{$\gm'$}
}
\]
Then $\gm'$ passes through $m$ and $n$, w.l.o.g. first through $m$
and then through $n$. But then $\gm'$ is not minimal because passing
through $a$, $c$ and $m$ (and possible further components between
$c$ and $m$) before passing through $n$ could be replaced by
just passing through $a$ and $b$ to arrive to $n$. By this contradiction
you are done here.

\case 
Both $b$, $b'$ are among $l$, $m$, $n$ and $o$. If $b$ and $b'$ are
neighbored, then $d(\gm)\le 3$.
\end{caselist}

Therefore, by this case distinction you are done unless at any
crossing of $\gm$ except the first and the last one $b$ and $b'$
participate as opposite components. In particular $b$ participates
as a neighboring component at any crossing of $\gm$ except possibly
the last one. But then one can see that $\gm$ must look like 
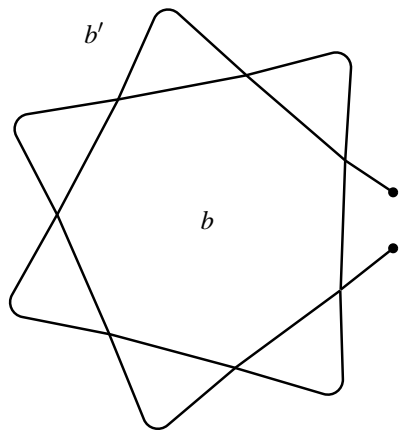
\[
\diag{1cm}{6}{6}{
  \pictranslate{3 3}{
    \picmultigraphics[rt]{6}{360 7 :}{
       \picline{2 23.5 polar}{2 75 polar}
       \picarcangle{2 23.5 polar}{3 50 polar}{2 75 polar}{0.2}
    }
    \picline{2 23.5 polar}{2.5 420 360 7 : - polar}
    \picline{2 -27.5 polar}{2.5 420 360 7 : - N polar}
    \picline{2 23.5 polar}{2 -27.5 polar}
    \ptat{2.5 420 360 7 : - polar}
    \ptat{2.5 420 360 7 : - N polar}
    \picputtext{0 0}{$b$}
    \picputtext{-1.5 2.5}{$b'$}
  }
}\,.
\]
To see this, start with 
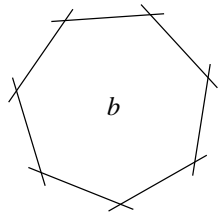
\[
\diag{5mm}{6}{6}{
  \pictranslate{3 3}{
    \picmultigraphics[rt]{7}{360 7 :}{
      \picline{2.8 10 polar}{2.8 75 polar}
    }
    \picputtext{0 0}{$b$}
  }
}
\]
and then observe that 
there is only one way to reconnect the outer arcs not creating
crossings (except possibly the last one) and having $b'$ as
specified, and moreover it works only if the number of crossings is odd.

But for such a curve $d(\gm)=0$. \qed

\section{Non-$Q$-maximal knots}

Now we are prepared to exhibit the Perko knot as non-$Q$-maximal.

\begin{theorem}
If $D$ is a prime diagram of a knot $K$ of $c(D)$ crossings with a
bridge of length $l=c(D)-k$, and $D$ has minimal crossing number among
all such diagrams for fixed $k$, then $l\le \max(3,k-3)$, hence
$c(D)\le k+\max(3,k-3)$.
\end{theorem}

From this we have the desired example:

\begin{exam}
If $10_{161}$ were $Q$-maximal, then we could pose $k=6$ in the
theorem and would obtain a 9 crossing diagram of the knot, which
does not exist. Hence $10_{161}$ is not $Q$-maximal.
\end{exam}

\proof[of theorem] This is basically lemma \reference{l2}.
Consider $\gm'$ to be the part of $D$ consisting of the maximal (length)
bridge and $\gm$ consisting of the rest of (the solid line of) $D$
with signs of all crossings ignored. Then the freedom to move the bridge
corresponds to the freedom to move $\gm'$. \qed

Clearly, for many phenomena Rolfsen's tables up to 10 crossings are very
limited. Verifying the list of non-alternating knots of at most 15
crossings provided by Thistlethwaite (see \cite{HTW}), I found 189 15 crossing
knots for which $\md Q\le 8$, and hence for which we would be done 
showing non-$Q$-maximality already with lemma \reference{l1} (or even
exercise \reference{x1}). The most striking examples are the knots
$15_{119574}$ and $15_{119873}$, where $\md Q=4$ (although for both
$\md_z F(a,z)=11$, the coefficients of the 7 highest powers of $z$
cancel when setting $a=1$).

There are several ways how the theorem can be modified.

\begin{theorem}
If $D$ is a diagram of a knot $K$ of $c(D)$ crossings with a
bridge of length $l=c(D)-k$, then $u(K)\le \br{k/2}$, where $u(K)$
denotes the unknotting number of $K$.
\end{theorem}

\proof By switching at most half of the crossings in $D$ not involved
in the maximal bridge, the remaining part $\gm$ of the plane curve 
(this time \em{with} signs of the crossings) can be layered, i.~e.,
any crossing is passed the first time as over- and then as under-%
crossing or vice versa. But reinstalling the bridge to a layered $\gm$
gives a layered, and hence unknotted, diagram. \qed

\begin{corr}\label{c2}
If $u(K)>\br{\md Q(K)/2}$, then $K$ is not $Q$-maximal. \qed
\end{corr}

Unfortunately, this corollary does not work to show non-$Q$-maximality
of Perko's knot. Verifying both hand-sides of the inequality (using that
the unknotting number of $10_{161}$ is 3, see \cite{pos,Kawamura,%
Tanaka}), we find that we just have equality. And that equality does
not suffice is seen, e.~g., from all 8 closed positive braid knots
in Rolfsen's tables (see \cite{Cromwell,Busk}) and more generally from
the $(2,n)$-torus knots for $n$ odd.

For knots of $>10$ crossings unknotting numbers are not tabulated
(anywhere I know of) and a general machinery does not exist to compute
them, hence when wanting to extend the search space for examples
applicable to corollary \reference{c2}, it makes sense to replace the
unknotting number by lower bounds for it, which can be computed
straightforwardly. I tried two such bounds. First we have the signature
$\sg$.

\begin{corr}
If $|\sg(K)|>\md Q(K)$, then $K$ is not $Q$-maximal. \qed
\end{corr}

Clearly, replacing $Q$ by lower bounds for it makes the condition more
and more restrictive. However, when checking the above mentioned
list of $189$ knots, I found that at least one of them satisfied
strict inequality. It is $15_{166028}$, where $\sg=8$ and $\md Q=7$.

\begin{figure}[htb]
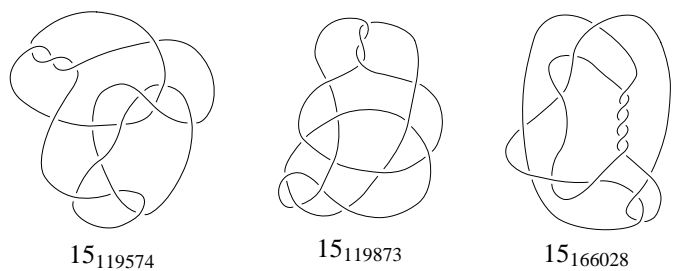

\[
\begin{array}{*3c}
\vis{15}{119574} &
\vis{15}{119873} &
\vis{15}{166028} \\
\end{array}
\]
\caption{\label{fignm}Three non-$Q$-maximal knots.}
\end{figure}

Another possibility is to minorate $u(K)$ by the bound coming from the
$Q$ polynomial itself.

\begin{corr}\label{c3}
If $2\log_{-3}Q(-1)>\md Q(K)$, then $K$ is not $Q$-maximal. \qed
\end{corr}

\begin{rem}
The negative logarithm base may disturb the reader because such logarithms
are usually not defined. But by work of Sakuma, Murakami, Nakanishi
(see Theorem 8.4.8 (2) of \cite{Kawauchi}) and Lickorish and
Millett \cite{LickMil} $Q(-1)$ is always a(n integral) power of $-3$
and this one it is referred to by this expression.
\end{rem}

The inequality in corollary \reference{c3} looks rather bizarre.
First, the inequality $u(K)\ge \log_{-3}Q(-1)$ is in general much less
sharp than the one with the signature and secondly, the inequality
in corollary \reference{c3} requires the coefficients of $Q$ to be of
an average magnitude which grows exponentially with $\md Q$. Thus,
non-surprisingly, my quest for applicable examples among the non-alternating
15 and 16 crossing knots ended with no success in this case.

\begin{question}
Is there a knot $K$ with $2\log_{-3}Q(-1)>\md Q(K)$?
\end{question}

I nevertheless gave the above inequality, because it is self-contained
w.r.t. $Q$ and would decide about non-$Q$-maximality
from $Q$ itself (without knowing anything else about the knot)
and hence is, in some sense, also beautiful.

\section{A question on plane curves}

The machinery of the dependence of $d(\gm)$ on $c(\gm)$ we developed
just as far as necessary for our knot theoretical context, but possibly
it is also interesting in its own right.

\begin{question}
Which is the best upper bound for $d(\gm)$ in terms of $c(\gm)$, i.~e.
a function $f:\,\bN\to\bN$ such that for any $\gm$ we have $f(c(\gm))\ge
d(\gm)$?
\end{question}

We proved that for \em{prime} $\gm$ we can choose $f(n):=\max(3,n-3)$.
Turning back to our example $\epsfs[6mm]{t-curve2}$, where $c=2$ and $d=1$,
and applying connected sums we find that we cannot choose $f(n)$ better
than $\br{n/2}$ (similar prime examples as
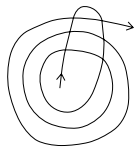
\[
\epsfs[17mm]{t-curve3}
\]
exist as well). But possibly this indeed is the best upper bound.
Unfortunately, proving it seems a matter of further tricky labour as in
\S\reference{S2}.

\end{document}